\documentclass[12pt]{article}
\usepackage{amssymb}
\newtheorem{dfn}{Definition}
\newtheorem{prop}{Proposition}
\newtheorem{thm}{Theorem}
\newtheorem{lem}{Lemma}
\newtheorem{cor}{Corollary}
\def\h{{\mathfrak{h}}}
\def\s{{\mathfrak{S}}}
\def\b{{\cal{B}}}

\title{A note on quantization operators on Nichols algebra model for 
Schubert calculus on Weyl groups}
\date{}
\author{Anatol N. Kirillov and  Toshiaki Maeno 
\vspace{8mm} \\ 
\small\it{ 
Dedicated to Kyoji Saito on the occasion of his sixtieth 
birthday }}
\begin{document}
\maketitle
\footnote{Both of the authors were 
supported by Grant-in-Aid for Scientific Research.}
\begin{abstract}
We give a description of the (small) quantum cohomology ring of the flag variety 
as a certain commutative subalgebra in the tensor product of the Nichols 
algebras. Our main result can be considered as a quantum analog of a result by 
Y. Bazlov. 
\end{abstract}
\section*{Introduction}
In this paper, we give a description of the (small) quantum cohomology rings of 
the flag varieties in terms of the braided differential calculus. Here, we give 
some remarks on the preceding works on this subject. 
In \cite{FK}, Fomin and one of the authors gave a combinatorial description of 
the Schubert calculus of the flag variety $Fl_n$ of type $A_{n-1}.$ 
They introduced a noncommutative quadratic algebra ${\cal E}_n$ determined 
by the root system, which contains 
the cohomology ring of the flag variety $Fl_n$ as a commutative subalgebra. 
One of remarkable 
properties of the algebra ${\cal E}_n$ is that it admits the quantum deformation, 
and the deformed algebra $\tilde{\cal E}_n$ also contains the quantum cohomology ring 
of the flag variety $Fl_n$ as its commutative subalgebra. 
A generalization of the algebras ${\cal E}_n$ and 
$\tilde{\cal E}_n$ was introduced by the authors in \cite{KM}. On the other hand, 
Fomin, Gelfand and Postnikov introduced the quantization operator on the polynomial 
ring to obtain the quantum deformation of the Schubert polynomials. Their approach 
was generalized for arbitrary root systems by Mar\'e \cite{Mar}. Our main idea is 
to lift their quantization operators onto the level of the Nichols algebras. 

The term ``Nichols algebra'' was introduced by 
Andruskiewitsch and Schneider \cite{AS}. 
The similar object was also discovered by Woronowicz \cite{Wo} and Majid \cite{Maj1} 
in the context of the braided differential calculus. 
The relationship between the quadratic algebra ${\cal E}_n$ 
and the Nichols algebra $\b(V_W)$ associated to a certain 
Yetter-Drinfeld module $V_W$ over the 
Weyl group $W$ was pointed out by Milinski and Schneider \cite{MS}. 
Majid \cite{Maj2} showed that it relates to a noncommutative differential structure 
on the permutation group $S_n.$ In fact, the higher order differential structure 
on $S_n$ gives a ``super-analogue'' of the algebra ${\cal E}_n.$ 
Recently, Bazlov \cite{Ba} showed that the Nichols algebra $\b(V_W)$ 
contains the coinvariant algebra 
${\bf S}_W$ of the finite Coxeter group $W.$ His method is based on 
the correspondence between braided derivations on $\b(V_W)$ and divided difference 
operators on the polynomial ring. Conjecturally, the algebra ${\cal E}_n$ is isomorphic 
to the Nichols algebra $\b(V_W)$ for $W=S_n.$ 
Our aim is to quantize his model for the coinvariant 
algebra in case $W$ is the Weyl group. 

Fix $B$ a Borel subgroup of a semisimple Lie group $G.$ 
Let $\h$ be the Cartan subalgebra in the Lie algebra of $G.$ 
We regard $\h$ as the 
reflection representation of the Weyl group $W.$ We have a set of positive roots 
$\Delta_+$ in the set of all roots $\Delta \subset \h^*.$ 
Denote by $\Sigma$ the set of 
simple roots. We need symbols 
$q^{\alpha^{\vee}}$ corresponding to the simple roots $\alpha$ as the parameters for 
the quantum deformation. Let $R$ be the polynomial ring 
${\bf C}[q^{\alpha^{\vee}} | \alpha \in \Sigma ].$ 
We also consider the algebra $\tilde{\b}(V)$ with a modified multiplication, 
see Section 1. 
Then our main result is: \smallskip \\ 
{\bf Theorem.} The algebra $(\b(V_W)\otimes \tilde{\b}(V_W))\otimes R$ contains 
the quantum cohomology ring of the flag variety $G/B$ as a commutative 
subalgebra. 
\medskip \\ 
{\bf Acknowledgements.} The authors would like to thank Yuri Bazlov for explaining 
his work and for fruitful comments. 
\section{Preliminaries} 
The aim of this paper is to describe the quantum cohomology ring of the flag 
variety in terms of the braided differential calculus. The Nichols algebra 
provides a suitable framework to consider the braided differential calculus. 
Let us recall some basic definitions. More detailed exposition can be found in 
\cite{Ba} and \cite{Maj1+}. 

Let $V$ be a finite dimensional ${\bf C}$-vector space equipped with a 
braiding $\Psi,$ i.e. a linear automorpism $\Psi : V \otimes V \rightarrow 
V \otimes V,$ subject to the braid relation 
\[ \Psi_{12}\Psi_{23}\Psi_{12}=\Psi_{23}\Psi_{12}\Psi_{23} \; \; \; \; {\rm on} 
\; \; V \otimes V \otimes V, \] 
where $\Psi_{ij}:V^{\otimes 3}\rightarrow V^{\otimes 3}$ stands for 
an automorphism obtained by applying $\Psi$ 
on the $i$-th and $j$-th components. 
The tensor algebra $T(V)$ of $V$ has a braided Hopf algebra structure with 
respect to the braiding induced by $\Psi.$ 
The coproduct $\triangle,$ the counit $\varepsilon$ and the antipode $S$ 
are defined by 
\[ \triangle (v) = v\otimes 1 + 1\otimes v, \; \varepsilon (v)=0, \; S(v)=-v, 
\; \; v \in V. \]
The dual space $V^*$ is a braided vector space with a braiding induced by $\Psi,$ 
and its tensor algebra $T(V^*)$ also has a structure of 
the braided Hopf algebra. 
The pairing $\langle \; , \; \rangle : V^* \times V \rightarrow {\bf C}, \; 
(\xi , x) \mapsto \xi(x)$ can be extended to the duality pairing 
$\langle \; , \; \rangle :T(V^*) \times T(V) \rightarrow {\bf C}$ 
so that the conditions 
\[ \langle \xi \eta , x \rangle = \langle \xi , x_{(2)} \rangle 
\langle \eta , x_{(1)} \rangle , \; 
\langle \xi , xy \rangle = \langle \xi_{(2)} ,x \rangle 
\langle \xi_{(1)} , y \rangle , \] 
\[ \langle 1 ,x \rangle = \varepsilon (x), \; \langle \xi , 1 \rangle = 
\varepsilon (\xi), \; \langle S(\xi),x \rangle = \langle \xi , S(x) \rangle \] 
are satisfied. Here, we use Sweedler's notation $\triangle (a)=a_{(1)}\otimes 
a_{(2)}.$ Let $I(V)$ be the kernel of the duality pairing. Then the 
Nichols-Woronowicz (or Nichols) algebra ${\cal B}(V)$ associated to $V$ is 
defined by ${\cal B}(V)= T(V)/I(V).$ One also has the dual algebra 
${\cal B}(V^*)$ as the quotient of $T(V^*)$ by the kernel of the duality 
pairing. It is known that the Nichols algebra ${\cal B}(V)$ constructed 
above coincides with the one characterized by the properties in the 
following definition. 
\begin{dfn}{\rm (Andruskiewitsch and Schneider \cite{AS})} 
The Nichols algebra ${\cal B}(V)$ associated to a braided vector space 
$V$ is a braided graded Hopf algebra satisfying the conditions: \\ 
{\rm (1)} ${\cal B}(V)^0 = {\bf C},$ \\ 
{\rm (2)} ${\cal B}(V)^1 = V = $ {\rm the set of primitive elements in} 
${\cal B}(V),$ \\ 
{\rm (3)} ${\cal B}(V)^1$ generates ${\cal B}(V)$ as an algebra. 
\end{dfn}
Note that each element in $V$ determines braided derivations acting on 
${\cal B}(V^*),$ some of which play a central role in the Nichols algebra 
model for the (quantum) Schubert calculus, see Definition 3. 

In the subsequent construction, we use a particular braided vector space 
called the Yetter-Drinfeld module. Let $G$ be a finite group and $V$ a 
finite dimensional $G$-module over ${\bf C}.$ 
\begin{dfn} 
The $G$-module $V$ is called the Yetter-Drinfeld 
module if $V$ has a $G$-grading, i.e. $V=\bigoplus_{g\in G} V_g,$ and the 
compatibility condition $gV_h=V_{ghg^{-1}}$ is satisfied. 
\end{dfn}
A significance of the 
Yetter-Drinfeld module is that it is braided naturally. The braiding $\Psi$ 
is given by $\Psi (x\otimes y)=gy\otimes x$ for $x\in V_g$ and $y\in V.$ 

Now let us proceed to our main ingredient. 
Consider the Nichols algebra $\b(V)$ associated to the Yetter-Drinfeld module 
$V=\bigoplus_{\alpha \in \Delta_+} {\bf C}[\alpha]$ over the Weyl group $W.$ The symbols 
$[\alpha]$ are subject to the condition $[-\alpha]=-[\alpha],$ and the $W$-action 
on $V$ is defined by $w.[\alpha]=[w(\alpha)].$ The $W$-degree of $[\alpha]$ is a 
reflection $s_{\alpha}\in W.$ The Yetter-Drinfeld module $V$ is a naturally braided 
vector space with a braiding $\psi_{V,V}.$ 
We can identify $\b(V)$ with its dual 
algebra $\b(V^*)$ via the $W$-invariant pairing $\langle [\alpha], [\beta] \rangle =
\delta_{\alpha,\beta}$ for $\alpha,\beta\in \Delta_+.$ 
Denote by $\tilde{\b}(V)$ the algebra $\b(V)$ with a modified 
multiplication $a*b=m(\psi^{-1}_{\b(V),\b(V)}(a\otimes b)),$ where $m$ is the multiplication map 
in the Nichols algebra $\b(V).$ 
\begin{dfn}
For each positive root $\alpha,$ the twisted derivation $\bar{D}_{\alpha}$ acting 
on $\b(V)$ from the left is defined by the rule 
\[ \bar{D}_{\alpha}([\beta])= \delta_{\alpha,\beta}, \; \; \; \beta \in \Delta_+, \] 
\[ (\dagger) \; \; \; \bar{D}_{\alpha}(xy)= \bar{D}_{\alpha}(x)y+
s_{\alpha}(x)\bar{D}_{\alpha}(y). \]  
\end{dfn}
The algebra $\tilde{\b}(V)$ acts on $\b(V^*)$ as an algebra generated by 
twisted derivations, and 
the twisted Leibniz rule $(\dagger)$ determines the algebra structure on 
$\b(V^*)\otimes \tilde{\b}(V):$ 
\[ (x\otimes [\alpha])\cdot (u\otimes v) = x\bar{D}_{\alpha}(u)\otimes v + xs_{\alpha}(u) 
\otimes [\alpha] *v .\] 
\begin{lem}
The representation of the algebra $\b(V^*)\otimes \tilde{\b}(V)$ on $\b(V^*)$ given by 
\[ ([\alpha_1] \cdots [\alpha_i]\otimes [\beta_1]*\cdots *[\beta_j] )(x) := 
[\alpha_1] \cdots [\alpha_i] \bar{D}_{\beta_1}\cdots \bar{D}_{\beta_j}(x) , \; \; \; 
x\in \b(V),  \] 
is faithful. 
\end{lem}
{\it Proof.} This follows from the non-degeneracy of the duality pairing between 
$\b(V^*)$ and $\b(V),$ {\it cf.} \cite{Ba}.  \rule{3mm}{3mm}  \medskip \\
Since the twisted derivations $\bar{D}_{\alpha}$ 
satisfy the Coxeter relations, one can define operators $\bar{D}_w$ for any element $w \in W$ 
by $\bar{D}_w=\bar{D}_{\alpha_1}\cdots \bar{D}_{\alpha_l}$ for a reduced decomposition 
$w=s_{\alpha_1}\cdots s_{\alpha_l}.$  
Let $R= {\bf C}[q^{\alpha^{\vee}} | \alpha \in \Delta_+ ] ,$ 
where the parameters $q^a$ satisfy the condition $q^{a+b}=q^aq^b.$ 
We denote by $\b_R(V)$ the scalar extension $R\otimes \b(V).$ 
Here, we define the quantization of the element $[\alpha] \in \b(V).$ 
Let $\tilde{\Delta}_+$ be the set of positive roots $\alpha$ satisfying the 
condition $l(s_{\alpha})=2{\rm ht}(\alpha^{\vee})-1,$ where the height ${\rm ht}(\alpha^{\vee})$ 
is defined by ${\rm ht}(\alpha^{\vee})=m_1+\cdots + m_n$ if 
$\alpha^{\vee}=m_1 \alpha_1^{\vee}+ \cdots + m_n \alpha_n^{\vee},$ $\alpha_i \in \Sigma.$ 
\begin{dfn} Let $(c_{\alpha})_{\alpha \in \Delta}$ be a set of nonzero constants 
with the condition $c_{\alpha}=c_{w\alpha},$ $w\in W.$ 
For each root $\alpha \in \Delta_+,$ we define an element 
$\widetilde{[\alpha]} \in \b_R(V^*)\otimes_R \tilde{\b}_R(V)$ 
by 
\[ \widetilde{[\alpha]} := \left\{ 
\begin{array}{cc}  
c_{\alpha}[\alpha] \otimes 1 +d_{\alpha}q^{\alpha^{\vee}} 
\otimes [\alpha_1] * \cdots * [\alpha_l], & {\rm if} \; \; \; 
\alpha \in \tilde{\Delta}_+, \\
c_{\alpha}[\alpha] \otimes 1, & {\rm otherwise,} 
\end{array} \right. \] 
where $\alpha_1,\ldots ,\alpha_l$ are simple roots appearing in a 
reduced decompositon $s_{\alpha}=s_{\alpha_1}\cdots s_{\alpha_l},$ and 
$d_{\alpha}=(c_{\alpha_1}\cdots c_{\alpha_l})^{-1}.$ 
We identify $\widetilde{[\alpha]}$ with an operator 
$c_{\alpha}[\alpha]  +d_{\alpha}q^{\alpha^{\vee}} \bar{D}_{s_{\alpha}}$ or 
a multiplication operator $c_{\alpha}[\alpha]$ 
acting on $\b_R(V^*)$ by Lemma 1. 
\end{dfn}
We define an $R$-linear map $\tilde{\mu} : \h_R \rightarrow V_R \otimes_R \b_R(V^*)$ 
in similar way to Bazlov \cite{Ba}, i.e., 
\[ \tilde{\mu}(x) =\sum_{\alpha \in \Delta_+}(x,\alpha)\widetilde{[\alpha]} . \] 
\begin{prop}
The subalgebra of $\b_R(V^*)\otimes_R \tilde{\b}_R(V)$ generated by ${\rm Im}(\tilde{\mu})$ is 
commutative. 
\end{prop}
{\it Proof.} We have to show $\tilde{\mu}(x)\tilde{\mu}(y)=
\tilde{\mu}(y)\tilde{\mu}(x)$ for arbitrary $x,y \in \h.$ 
The left hand side is expanded as 
\[ (*) \; \; \; \; \sum_{\alpha, \beta \in \Delta_+}(x,\alpha)(y,\beta) 
c_{\alpha}c_{\beta}[\alpha][\beta] \]
\[ + \sum_{\alpha \in \tilde{\Delta}_+, \beta \in \Delta_+}
(x,\alpha)(y,\beta) d_{\alpha}c_{\beta}
q^{\alpha^{\vee}}\bar{D}_{s_{\alpha}}\cdot [\beta] + 
\sum_{\alpha \in \Delta_+, \beta \in \tilde{\Delta}_+} (x,\alpha)(y,\beta)
c_{\alpha}d_{\beta}q^{\beta^{\vee}}[\alpha] \cdot \bar{D}_{s_{\beta}} \] 
\[ +\sum_{\alpha \in \tilde{\Delta}_+, \beta \in \tilde{\Delta}_+} (x,\alpha)(y,\beta)
d_{\alpha}d_{\beta}q^{\alpha^{\vee}+\beta^{\vee}}\bar{D}_{s_{\alpha}}\bar{D}_{s_{\beta}} . \]
We have already known the commutativity of the classical part (\cite{Ba}, \cite{KM}), 
so we can ignore the first 
term.  We also have 
\[ \bar{D}_{s_{\alpha}}\bar{D}_{s_{\beta}}=\left\{ 
\begin{array}{cc} 
\bar{D}_{s_{\alpha}s_{\beta}} & {\rm if} \; \; l(s_{\alpha}s_{\beta})=l(s_{\alpha})+l(s_{\beta}), \\ 
0 & {\rm otherwise,} 
\end{array} \right. \]
and 
\[ \bar{D}_{s_{\alpha}}\cdot [\beta] - s_{\alpha}([\beta])\bar{D}_{s_{\alpha}}=\left\{ 
\begin{array}{cc}
\bar{D}_{s_{\alpha}s_{\beta}} & {\rm if} \; \; l(s_{\alpha}s_{\beta})=l(s_{\alpha})-1, \\ 
0 & {\rm otherwise.} 
\end{array} \right. \] 
Let 
\[ A=\{ (\alpha,\beta) \in \tilde{\Delta}_+ \times \Delta_+ | 
\; l(s_{\alpha}s_{\beta})=l(s_{\alpha})-1 \}  \]  
and 
\[ B=\{ (\alpha,\beta)\in \tilde{\Delta}_+^2 | \; l(s_{\alpha}s_{\beta})=
l(s_{\alpha})+l(s_{\beta}) \} . \] 
Then, we have 
\[ \sum_{\alpha \in \tilde{\Delta}_+, \beta \in \Delta_+}
(x,\alpha)(y,\beta) d_{\alpha}c_{\beta}
q^{\alpha^{\vee}}\bar{D}_{s_{\alpha}}\cdot [\beta] + 
\sum_{\alpha \in \Delta_+, \beta \in \tilde{\Delta}_+} (x,\alpha)(y,\beta)
c_{\alpha}d_{\beta}q^{\beta^{\vee}}[\alpha] \cdot \bar{D}_{s_{\beta}} \] 
\[= \sum_{\alpha \in \Delta_+, \beta \in \tilde{\Delta}_+} 
c_{\alpha}d_{\beta}\left((x,\alpha)(y,\beta)+(x,\beta)(y,\alpha)-2
(\alpha,\beta)(x,\beta)(y,\beta) \right)
q^{\beta^{\vee}}[\alpha]\cdot \bar{D}_{s_{\beta}} \] 
\[ + \sum_{(\alpha, \beta) \in A} 
d_{\alpha}c_{\beta}(x,\alpha)(y,\beta)q^{\alpha^{\vee}}\bar{D}_{s_{\alpha}s_{\beta}} , \] 
and 
\[ \sum_{\alpha, \beta \in \tilde{\Delta}_+} d_{\alpha}d_{\beta}(x,\alpha)(y,\beta) 
q^{\alpha^{\vee}+\beta^{\vee}}\bar{D}_{s_{\alpha}}\bar{D}_{s_{\beta}}
= \sum_{(\alpha, \beta) \in B} d_{\alpha}d_{\beta}(x,\alpha)(y,\beta) 
q^{\alpha^{\vee}+\beta^{\vee}}\bar{D}_{s_{\alpha}s_{\beta}} . \]
For each element $(\alpha,\beta)\in A$ with $\alpha \not= \beta,$ we can find 
an element $(\gamma,\delta)\in B$ such that $\alpha^{\vee}=\gamma^{\vee}+\delta^{\vee}$ 
and $s_{\alpha}s_{\beta}=s_{\gamma}s_{\delta}$ 
from the argument in \cite[Section 3]{Mar}. 
This correspondence gives a bijection between the sets 
$A'=A\setminus \{ (\alpha,\beta) | \alpha = \beta \}$ and $B'=B\setminus 
\{ (\gamma,\delta) | s_{\gamma}s_{\delta}= s_{\delta}s_{\gamma} \},$ 
and $(x,\alpha)(y,\beta)+(x,\gamma)(y,\delta)$ is symmetric in $x$ and $y$ under 
the correspondence between $(\alpha,\beta)\in A'$ and $(\gamma,\delta)\in B'.$ 
Hence, $(*)$ is symmetric in $x$ and $y.$  \rule{3mm}{3mm} \medskip \\ 
{\bf Remark.} We can use the opposite algebra $\b(V)^{op}$ and the twisted 
derivation $\overleftarrow{D}_{\alpha}$ acting from the right, instead of 
$\tilde{\b}(V)$ and $\bar{D}_{\alpha}.$ The algebra $\b(V)^{op}$ is 
the opposite algebra of $\b(V),$ whose 
multiplication $\star$ is obtained by reversing the order of the multiplication in $\b(V),$ 
i.e., 
\[ a_1 \star \cdots \star a_m = a_m \cdots a_1. \] 
The twisted derivation $\overleftarrow{D}_{\alpha},$ $\alpha \in \Delta_+,$ is determined 
by the conditions: 
\[ [\beta]\overleftarrow{D}_{\alpha} = \delta_{\alpha,\beta} , \; \; \; \beta 
\in \Delta_+, \] 
\[ (fg)\overleftarrow{D}_{\alpha} =f(g\overleftarrow{D}_{\alpha})+
(f\overleftarrow{D}_{\alpha})s_{\alpha}(g). \] 
Then, the algebra $\b(V^*) \otimes \b(V)^{op}$ faithfully acts on the 
algebra $\b(V^*)$ from the left via $1\otimes [\alpha] \mapsto \overleftarrow{D}_{\alpha}$ 
and $[\beta]\otimes 1\mapsto$ (left multiplication by $[\beta]$). 
We can also define the quantized element 
$\widetilde{[\alpha]}$ as an element in $\b_R(V^*) \otimes_R \b_R(V)^{op}$ in a similar 
way to Definition 4: 
\[ \widetilde{[\alpha]} := \left\{ 
\begin{array}{cc} 
c_{\alpha}[\alpha] \otimes 1 +d_{\alpha}q^{\alpha^{\vee}}\otimes 
[\alpha_1]\star \cdots \star [\alpha_l] , 
& {\rm if} \; \; \; 
\alpha \in \tilde{\Delta}_+, \\
c_{\alpha}[\alpha] \otimes 1, & {\rm otherwise.} 
\end{array} \right. \]  
The arguments in this section work well for this definition, in particular, 
the subalgebra generated by ${\rm Im}(\tilde{\mu})$ is again commutative. 
This construction of 
the quantized elements $\widetilde{[\alpha]}$ by using 
$\b_R(V)^{op}$ and the twisted derivations from the right was suggested by Bazlov. 
\section{Main result}
Now we can extend $\tilde{\mu}$ as an $R$-algebra homomorphism 
${\rm Sym_R}(\h_R) \rightarrow \b_R(V^*)\otimes _R \tilde{\b}_R(V).$  
Let $\mu : {\rm Sym}_R(\h_R)\rightarrow \b_R(V^*)$ be the scalar extension of 
the homomorphism introduced in \cite{Ba}, i.e., 
\[ \mu(x)= \sum_{\alpha \in \Delta_+} c_{\alpha}(x,\alpha)[\alpha]. \]
The Demazure operator $\partial_{\alpha},$ $\alpha \in \Delta_+,$ acting on 
the polynomial ring ${\rm Sym}(\h)$ is defined by $\partial_{\alpha}(f)=
(f-s_{\alpha}(f))/\alpha.$ For each element $w\in W,$ the operator 
$\partial_w$ can be defined as $\partial_w=\partial_{\alpha_1}\cdots 
\partial_{\alpha_l}$ for a reduced decomposition $w=s_{\alpha_1}\cdots 
s_{\alpha_l},$ $\alpha_1,\ldots ,\alpha_l \in \Sigma.$ This is well-defined 
since the Demazure operators satisfy $\partial_{\alpha}^2=0$ and the Coxeter relations.
\begin{lem} {\rm (\cite{Ba})} For $f \in {\rm Sym}(\h),$ we have 
\[ \bar{D}_{\alpha}\mu (f) = c_{\alpha} \mu (\partial_{\alpha}f). \] 
\end{lem}
\begin{prop}
Let $I_i^q,$ $1\leq i \leq n={\rm rk}\h,$ be the quantum fundamental $W$-invariants 
given by {\rm \cite{GK}} and {\rm \cite{Ki}.} 
Then, $\tilde{\mu}(I_i^q)\mu(f)=0,$ $\forall f \in {\rm Sym}_R(\h_R).$  
\end{prop}
{\it Proof.} 
For each simple root $\alpha \in \Sigma,$ we define 
\[ \eta_{\alpha}:= \sum_{\gamma \in \Delta_+}\langle \omega_{\alpha}, \gamma^{\vee} 
\rangle \widetilde{[\gamma]}
=\sum_{\gamma \in \Delta_+}\langle \omega_{\alpha}, \gamma^{\vee} \rangle
c_{\gamma}[\gamma]+\sum_{\gamma \in \tilde{\Delta}_+}\langle \omega_{\alpha}, \gamma^{\vee} \rangle
d_{\gamma}q^{\gamma^\vee}\bar{D}_{s_{\gamma}}, \]
where $\omega_{\alpha}$ is a fundamental dominant weight corresponding to $\alpha.$ 
Then, Lemma 2 shows that 
\[ \eta_{\alpha} \mu(f) = \mu (Y_{\alpha} f), \]
where 
\[ Y_{\alpha}= \omega_{\alpha} + 
\sum_{\gamma \in \tilde{\Delta}_+}\langle \omega_{\alpha}, \gamma^{\vee} \rangle
q^{\gamma^{\vee}}\partial_{s_{\gamma}} . \] 
Hence, $\tilde{\mu}(\varphi) \mu(f) = \mu( \varphi((Y_{\alpha})_{\alpha})(f))$ 
for any polynomial $\varphi \in {\rm Sym}_R (\h_R).$ From the quantum Pieri or Chevalley 
formula (\cite{FGP}, \cite{FW}, \cite{P}), we 
have $\tilde{\mu}(I_i^q)(1)=0.$ For any $f\in {\rm Sym}_R(\h_R),$ there exists 
a polynomial $\tilde{f}\in {\rm Sym}_R(\h_R)$ such that 
$\tilde{f}((Y_{\alpha})_{\alpha})(1)=f.$ 
Then, we have 
\[ \tilde{\mu}(I_i^q)\mu(f)=\tilde{\mu}(I_i^q)\tilde{\mu}(\tilde{f})(1)=
\tilde{\mu}(\tilde{f})\tilde{\mu}(I_i^q)(1)=0. \; \; \; \rule{3mm}{3mm} \] 
\begin{thm}
${\rm Im}(\tilde{\mu})$ generates a subalgebra in $\b_R(V^*)\otimes_R \tilde{\b}_R(V)$ 
isomorphic to the quantum cohomology ring of the corresponding flag variety $G/B.$ 
\end{thm}
{\it Proof.} We assign the degree 1 to the elements $[\alpha]$ and $-1$ to 
$\bar{D}_{\alpha}.$ 
Define the filter $F_{\bullet}$ on the algebra ${\rm Im}(\tilde{\mu})$ by 
$F_i({\rm Im}(\tilde{\mu}))= \{ x | \deg(x) \leq i \} .$ Then, 
$Gr_F({\rm Im}{\tilde{\mu}}) \cong {\rm Im}(\mu).$ 
The faithfulness of the representation of the subalgebra 
${\rm Im}(\mu)$ in $\b_R(V)$ on itself implies that of the representation of the algebra 
generated by ${\rm Im}(\tilde{\mu})$ on ${\rm Im}(\mu).$ Hence, we have 
$\tilde{\mu}(I_i^q)=0$ from Proposition 2. Since $Gr_F({\rm Im}{\tilde{\mu}}) 
\cong {\rm Im}(\mu),$ we conclude that 
${\rm Im}(\tilde{\mu}) \cong {\rm Sym}_R(\h_R)/(I_1^q,\ldots,I_n^q).$ \rule{3mm}{3mm} 
\begin{cor}
{\rm (1)} In the case of root systems of type $A_n,$ denote by $\s_w$ and $\s_w^q$ the 
Schubert polynomial and its quantization corresponding to $w\in S_{n+1}.$ 
Then, $\tilde{\mu}(\s_w^q)(1)= \mu(\s_w).$ \\
{\rm (2)} For general crystallographic root systems, let $X_w$ and $X_w^q$ be the 
Bernstein-Gelfand-Gelfand polynomial {\rm (\cite{BGG})} and its quantization 
coresponding to $w\in W$ {\rm (\cite{KM},\cite{Mar}).} 
Then, $\tilde{\mu}(X_w^q)(1)= \mu(X_w).$ 
\end{cor}
{\bf Remark.} In $A_n$-cases, the operators $\eta_{\alpha}$ induce 
the operators on the algebra ${\rm Sym}_R (\h_R)$ introduced by Fomin, 
Gelfand and Postnikov \cite{FGP}. 
For other cases, they 
induce Mar\'e's operators \cite{Mar}. The above corollary is a restatement of their 
results and \cite[Proposition 8.1]{KM}. 
\begin{prop} The identity 
\[ \widetilde{[\alpha]}^2= \left\{ 
\begin{array}{cc}
c_{\alpha}d_{\alpha}q^{\alpha^{\vee}}, & \rm{if} \; \; \; \alpha : \; \rm{simple,} \\
0, & \rm{otherwise} 
\end{array}
\right. \]
holds in $\b_R(V^*)\otimes_R \tilde{\b}_R(V).$ 
\end{prop}
{\it Proof.} This follows from $[\alpha]^2=0,$ $\bar{D}_{s_{\alpha}}^2=0$ and 
\[ \bar{D}_{s_{\alpha}} \cdot [\alpha]  = \left\{ 
\begin{array}{cc}
1-[\alpha]\bar{D}_{s_{\alpha}}, & \rm{if} \; \; \; \alpha : \; \rm{simple,} \\
-[\alpha]\bar{D}_{s_{\alpha}}, & \rm{otherwise.} 
\end{array}
\right. \] 
{\bf Example.} In $B_n$-case, the algebra $\b(V)$ is generated by the symbols 
$[i,j],$ $\overline{[i,j]}$ and $[i]$ with $1\leq i,j \leq n$ and $i\not= j.$ 
After normalizing $c_{\alpha}=1$ for all $\alpha \in \Delta,$ 
the quantized operators are given by 
\[ \widetilde{[i,j]} = [i,j] + Q_{ij}\bar{D}_{(ij)}, \; \; \; (i<j) , \] 
\[ \widetilde{\overline{[i,j]}} = \overline{[i,j]} + 
Q_{\overline{ij}}\bar{D}_{\overline{(ij)}} , \] 
\[ \widetilde{[i]} = [i], \; \; \; (i<n) , \] 
\[ \widetilde{[n]} = [n] + Q_{n}\bar{D}_{(n)} , \]
where $Q_{ij}=q_iq_j^{-1}$ ($i<j$), $Q_{\overline{ij}}=q_iq_j$ and 
$Q_n=q_n^2$ are elements in the Laurent polynomial ring 
${\bf C}[q_1^{\pm 1},\ldots , q_n^{\pm 1}].$ We put 
$\widetilde{[j,i]}=-\widetilde{[i,j]}.$ 
We can check that $\widetilde{[i,j]},$ 
$\widetilde{\overline{[i,j]}}$ and $\widetilde{[i]}$ satisfy the relations 
of the quantum $B_n$-bracket algebra introduced by the authors \cite{KM}: \\
(1) $\widetilde{[i,i+1]}^2=Q_{i \; i+1},$ $\widetilde{[n]}^2=Q_n,$ \\ 
\hspace*{5.2mm} $\widetilde{[i,j]}^2=0$, if $| i-j | \not= 1$;
$\widetilde{[i]}^2=0$, if $i<n$; $\widetilde{\overline{[i,j]}}^2=0,$ 
if $i\not= j,$ \smallskip \\ 
(2) $\widetilde{[i,j]}\widetilde{[k,l]}=
\widetilde{[k,l]}\widetilde{[i,j]},$ 
$\widetilde{\overline{[i,j]}}\widetilde{[k,l]}=
\widetilde{[k,l]}\widetilde{\overline{[i,j]}},$ 
$\widetilde{\overline{[i,j]}}\widetilde{\overline{[k,l]}}=
\widetilde{\overline{[k,l]}}\widetilde{\overline{[i,j]}},$ \\ 
\hspace*{5.2mm} if $\{ i,j\} \cap \{ k,l\} =\o $, 
\smallskip \\ 
(3) $\widetilde{[i]}\widetilde{[j]}=
\widetilde{[j]}\widetilde{[i]},$ 
$\widetilde{[i,j]}\widetilde{\overline{[i,j]}}=
\widetilde{\overline{[i,j]}}\widetilde{[i,j]},$ 
$\widetilde{[i,j]}\widetilde{[k]}=\widetilde{[k]}\widetilde{[i,j]}$, 
$\widetilde{\overline{[i,j]}}\widetilde{[k]}=
\widetilde{[k]}\widetilde{\overline{[i,j]}}$, 
if $k\not= i,j$, \smallskip \\ 
(4) $\widetilde{[i,j]}\widetilde{[j,k]}+
\widetilde{[j,k]}\widetilde{[k,i]}+
\widetilde{[k,i]}\widetilde{[i,j]}=0,$ \\ 
\hspace*{5.2mm} $\widetilde{\overline{[i,k]}}\widetilde{[i,j]}+
\widetilde{[j,i]}\widetilde{\overline{[j,k]}}+
\widetilde{\overline{[k,j]}}\widetilde{\overline{[i,k]}}=0,$ \\ 
\hspace*{5.2mm} $\widetilde{[i,j]}\widetilde{[i]}+
\widetilde{[j]}\widetilde{[j,i]}+\widetilde{[i]}\widetilde{\overline{[i,j]}}+
\widetilde{\overline{[i,j]}}\widetilde{[j]}=0,$ \\
\hspace*{5.2mm} if all $i,$ $j$ and $k$ are distinct,
\smallskip \\  
(5) $\widetilde{[i,j]}\widetilde{[i]}\widetilde{\overline{[i,j]}}\widetilde{[i]}+
\widetilde{\overline{[i,j]}}\widetilde{[i]}\widetilde{[i,j]}\widetilde{[i]}+
\widetilde{[i]}\widetilde{[i,j]}\widetilde{[i]}\widetilde{\overline{[i,j]}}+
\widetilde{[i]}\widetilde{\overline{[i,j]}}\widetilde{[i]}\widetilde{[i,j]}=0,$ 
if $i<j.$ 
\medskip \\ 
{\bf Remark.} As in the remark at the end of Section 1, we also have another 
construction of the quantized elements by using $\b(V)^{op}$ and $\overleftarrow{D}_{\alpha}.$ 
Since 
\[ \mu (f)\overleftarrow{D}_{\alpha}=c_{\alpha}\mu(\partial_{\alpha}f) \] 
is also correct, we can show that the algebra $\b_R(V^*)\otimes_R \b_R(V)^{op}$ 
contains the quantum cohomology ring of $G/B$ as a commutative subalgebra. 

Research Institute for Mathematical Sciences, \\
Kyoto University, \\ 
Sakyo-ku, Kyoto 606-8502, Japan \\ 
e-mail: kirillov@kurims.kyoto-u.ac.jp 
\medskip \\
Department of Mathematics, \\
Kyoto University, \\ 
Sakyo-ku, Kyoto 606-8502, Japan \\ 
e-mail: maeno@math.kyoto-u.ac.jp

\begin{thebibliography}{99}
\bibitem{AS} N. Andruskiewitsch and H.-J. Schneider, {\it Pointed Hopf algebras,} 
New directions in Hopf algebras, Math. Sci. Res. Inst. Publ. {\bf 43}, 
Cambridge Univ. Press, Cambridge, 2002, 1-68. 
\bibitem{Ba} Y. Bazlov, {\it Nichols-Woronowicz algebra model for Schubert calculus on 
Coxeter groups,} preprint, math.QA/0409206. 
\bibitem{BGG} I. N. Bernstein, I. M. Gelfand and S. I. Gelfand, {\it Schubert 
cells and 
cohomology of the spaces $G/P,$} Russian Math. Surveys {\bf 28} (1973), 1-26.
\bibitem{FGP} S. Fomin, S. Gelfand and A. Postnikov, {\it Quantum Schubert polynomials,} 
J. Amer. Math. Soc. {\bf 10} (1997), 565-596. 
\bibitem{FK} S. Fomin and A. N. Kirillov, {\it Quadratic algebras, 
Dunkl elements and Schubert calculus,} Advances in 
Geometry, (J.-L. Brylinski, R. Brylinski, V. Nistor, B. Tsygan 
and P. Xu, eds.) Progress in Math. {\bf 172} Birkh\"auser, 
1995, 147-182. 
\bibitem{FW} W. Fulton and C. Woodward, {\it On the quantum 
product of Schubert classes,} J. Alg. Geom. {\bf 13} (2004), 641-661.
\bibitem{GK} A. Givental and B. Kim, {\it Quantum cohomology of flag 
manifolds and Toda lattices,} 
Commun. Math. Phys. {\bf 168} (1995), 609-641. 
\bibitem{Ki} B. Kim, {\it Quantum cohomology of flag 
manifolds $G/B$ and quantum Toda lattices,} 
Annals of Math. {\bf 149} (1999), 129-148. 
\bibitem{KM} A. N. Kirillov and T. Maeno, {\it Noncommutative algebras related with 
Schubert calculus on Coxeter groups,} European J. of Combin. {\bf 25} (2004), 1301-1325. 
\bibitem{Maj1} S. Majid, {\it Free braided differential calculus, braided binomial 
theorem, and the braided exponential map,} J. of Math. Phys. {\bf 34} (1993), 4843-4856. 
\bibitem{Maj1+} S. Majid, {\it Foundations of quantum group theory,} Cambridge University 
Press, 1995. 
\bibitem{Maj2} S. Majid, {\it Noncommutative differentials and 
Yang-Mills on permutation groups $S_N,$} 
Lect. Notes Pure Appl. Math. {\bf 239} (2004), 189-214.
\bibitem{Mar} A.-L. Mar\'e, {\it The combinatorial quantum 
cohomology ring of $G/B,$} preprint, math.CO/0301257.
\bibitem{MS} A. Milinski and H.-J. Schneider, {\it Pointed 
indecomposable Hopf algebras over Coxeter groups,} Contemp. 
Math. {\bf 267} (2000), 215-236. 
\bibitem{P} D. Peterson, {\it Quantum cohomology 
of flag varieties, Lectures given at MIT,} 1997. 
\bibitem{Wo} S. L. Woronowicz, {\it Differential calculus 
on compact matrix pseudogroups (quantum groups),} Commun. Math. 
Phys. {\bf 122} (1989), 125-170. 
\end{thebibliography}
\end{document}